\documentclass[letterpaper, 10 pt, conference]{ieeeconf}  

\IEEEoverridecommandlockouts                              
\usepackage{amssymb}
\usepackage{graphicx}
\graphicspath{{./figs/}}
\usepackage{float}
\usepackage{color}
\usepackage{mathrsfs} 
\usepackage{mathtools}
\usepackage{caption}
\usepackage{subcaption}
\usepackage{cite}
\usepackage{bm}
\usepackage{hyperref}
\usepackage{cleveref}

\usepackage{amsthm}

\theoremstyle{plain}
\newtheorem{theorem}{Theorem}

\newtheorem{remark}[theorem]{Remark}
\newtheorem{proposition}[theorem]{Proposition}

\theoremstyle{definition}

\newtheorem{assumption}[theorem]{Assumption}

\newcommand{\R}{\mathbb{R}}
\newcommand{\N}{\mathbb{N}}

\title{\LARGE \bf Model Predictive Control of Spreading Processes\\ via Sparse Resource Allocation}

\author{Ruigang Wang, Armaghan Zafar and Ian R. Manchester
	\thanks{This work was supported by the Air Force Office of Scientific Research Grant FA2386-19-1-4076 and the NSW Defence Innovation Network.}
	\thanks{The authors are with the Australian Centre for Field Robotics, and Sydney Institute for Robotics and Intelligent Systems, The University of Sydney, Sydney, NSW 2006, Australia
		(e-mail: {\tt\small  ruigang.wang, armaghan.zafar. ian.manchester@sydney.edu.au}).}%
}

\begin{document}

\maketitle
\pagestyle{empty}

\begin{abstract}
In this paper, we propose a model predictive control (MPC) method for real-time intervention of spreading processes, such as epidemics and wildfire, over large-scale networks. The goal is to allocate budgeted resources each time step to minimize the risk of an undetected outbreak, i.e. the product of the probability of an outbreak and the impact of that outbreak. By using dynamic programming relaxation, the MPC controller is reformulated as a convex optimization problem, in particular an exponential cone programming. We also provide sufficient conditions for the closed-loop risks to asymptotically decrease and a method to estimate the upper bound of when the risk will monotonically decrease. Numerical results are provided for a wildfire example.
\end{abstract}

\section{Introduction}

Modeling, analysis and control of epidemics and wildfires have been studied intensely due to their significant social and economic impacts, see the recent surveys \cite{nowzari2016analysis,pare2020modeling,zino2021analysis,karafyllidis1997model}. These events are usually modeled as spreading process in which an initial localized outbreak spreads rapidly to neighboring nodes in a network, and causes harm or risk of harm to the nodes it touches. The task of intervention is to reduce the risks of an outbreak by allocating resources (e.g. vaccines, social distance) to nodes/links. 

There are several challenges for resource allocation of large-scale spreading processes. Firstly, the available resources at any given time are often limited. Secondly, some types of resources (e.g. water bombing attacks on a wildfire) can be difficult or impossible to widely distribute across the network. Moreover, the large-scale nature of the spreading processes, e.g. global transport network for epidemics or large geographic areas for bushfires, imposes computational challenges in terms of real-time computation of responses. Therefore, it is important to develop scalable resource allocation methods resulting in sparse solutions, subject to meaningful budget constraints. 

A natural tool for developing such intervention strategies is optimal control with the objective of minimizing the risk of an outbreak, which is defined as the expected accumulated impact of the event over the time horizon. Typical models for spreading processes are the Susceptible-Infected-Susceptible (SIS) model and the Susceptible-Infected-Removed (SIR) model \cite{bailey1975mathematical}. These models are stochastic which can further increase the complexity of optimal control problem. By using mean-field approximation \cite{van2015accuracy}, one can obtain approximate models in the form of nonlinear deterministic ordinary differential equations (ODEs). It is proven \cite{van2008virus} that solutions of these deterministic ODE systems are upper bounds of the expected value of solutions of the corresponding stochastic models. Thus, they are usually the object of study in control design. 

Various methods have been proposed for resource allocation subject to budget constraints. Convex optimization techniques such as semidefinite programming (SDP) \cite{preciado2013optimal} and geometric programming (GP) \cite{preciado2014optimal,nowzari2015optimal} have been applied. In \cite{somers2020sparse}, a novel logarithmic resource model was proposed for encouraging sparsity solutions for the GP-based approaches. Further explorations in surveillance study of spreading processes can be found in \cite{somers2022minimizing}. However, all these methods offer static intervention strategies only. That is, the resource allocation is determined by a single-step optimization problem based on the current state instead of allocated over time. 

Dynamic resource allocation, in particular using Model Predictive Control (MPC), has been studied recently for epidemic intervention. In \cite{selley2015dynamic}, nonlinear MPC is applied to reduce the epidemic spread. The works \cite{kohler2018dynamic,kohler2021robust} investigate MPC approaches to control study of the COVID-19 outbreak in Germany. In \cite{watkins2019robust}, robust economic MPC with rigorous convergence guarantees was developed for the containment of stochastic epidemic processes. However, these approaches do not provide sparse resource allocation.

In \cite{somers2022multi}, a multi-stage sparse resource allocation method was developed by extending the single-stage approach \cite{somers2020sparse,somers2022minimizing}. Different from the MPC approaches \cite{selley2015dynamic,kohler2018dynamic,watkins2019robust,kohler2021robust}, this method allocates resources by minimizing a risk upper bound, which is a linear value (cost-to-go) function constructed by dynamic programming. The resulting optimal control problem (OCP) can be reformulated as a convex optimization problem, in particular an exponential cone program (ECP), which can be efficiently solved by the state-of-art tools (e.g. MOSEK \cite{mosek}). 

In this paper we, we build on \cite{somers2022multi} by providing a theoretical analysis of the closed-loop behavior when the this multi-stage resource allocation approach is applied recursively, in an MPC framework. In particular, we  prove recursive feasibility and provide sufficient conditions for the closed-loop risk bound to asymptotically decrease. Moreover, we formulate a convex optimization problem to estimate the time at which the risk will begin to monotonically decrease. We also provide results of the numerical experiments performed by implementing the proposed MPC approach on a wildfire example.

The remainder of this paper is structured as follows. The spreading process model and control setup are presented in Section~\ref{sec:problem}. Section~\ref{sec:dp} gives the proposed approach and some theoretical results with proofs found in Section~\ref{sec:proof}. Finally, an illustrative example of a wildfire is presented in Section~\ref{sec:example}.

\textbf{Notation.} Throughout the paper, $\geq$ is an element-wise operator, e.g., the matrices $A,B$ satisfying $A\geq B$ implies $A_{ij}\geq B_{ij}$ for all $(i,j)$. And $A\gneqq  B$ means $A_{ij}\geq B_{ij}$ for all $(i,j)$ and there exist some pairs $(i,j)$ such that $A_{ij}>B_{ij}$. We call $A$ a nonnegative matrix if $A\geq 0$.

\section{Problem Formulation}\label{sec:problem}
\subsection{Spreading Process over Networks}

We consider an SIS spreading process on a graph $\mathcal{G}(\mathcal{V},\mathcal{E})$ with node set $\mathcal{V}=\{1,2,\ldots,n\}$ and edge set $\mathcal{E}=\{(i,j):i\neq j\in\mathcal{V}\}$, where each node $i$ has a state $X_i(t)\in\{0,1\}$ indicating that it is infectious ($X_i(t)=1$) or susceptible ($X_i(t)=0$). In this process, an infected node $i$ recovers with probability $\delta_i(t) \Delta t$ to $X_i(t+\Delta t)=0$ and an susceptible node $i$ can be infected (i.e. $X_i(t+\Delta t)=1$) by its neighbor node $j$ with probability $\beta_{ij}(t)\Delta t$, where $\Delta t$ is a small time interval. By using mean-field approximation \cite{van2015accuracy}, we obtain $n$ coupled deterministic differential equations:
\begin{equation}\label{eq:sys-ct}
    \dot{x}_i(t)=(1-x_i(t))\sum_{j=1}^{n}\beta_{ij}x_j(t)-\delta_i x_i(t),
\end{equation}
where $x_i(t):=P(X_i(t)=1)$ is the probability of a node $i$ being infected at time $t$. 

By applying Euler forward method to \eqref{eq:sys-ct}, we obtain a discrete-time model in compact form of
\begin{equation}\label{eq:sys}
    x(k+1)=\bigl[A-B(x(k))\bigr]x(k),
\end{equation}
where $x(k)=\bigl[x_1(kh),\ldots,x_n(kh)\bigr]^\top$ is the state at time step $k\in\N$, and the matrices $A,B$ are given by
\begin{equation}\label{eq:AB}
    \begin{split}
        A_{ij}&=\begin{cases}
    1-h\delta_i, \; i=j, \\
    h\beta_{ij}, \; (i,j)\in \mathcal{E} ,
    \end{cases}
    B_{ij}(x)=
    hx_i\beta_{ij}.
    \end{split}
\end{equation}
We assume that the time interval $h>0$ satisfies $ h\delta_j \leq 1$ and $h\sum_{i=1}^{n}\beta_{ij}<1$ for all $j\in \mathcal{V}$. Note that \eqref{eq:sys} is well-defined, i.e., $x(0)\in [0,1]^n \Rightarrow x(k) \in [0,1]^n$ for all $k$.

\subsection{Risk Model}
We define the risk of an outbreak $x$ as the impact over the time horizon, i.e., 
\begin{equation}\label{eq:risk}
    R(x)=\sum_{\ell=0}^{\infty} \alpha^{k}C x(k),
\end{equation}
where $x(k)$ is state of \eqref{eq:sys} with initial condition $x(0)=x$. Here $C=\bigl[c_1,\ldots,c_n\bigr]>0$ is a row vector defining the cost associate with each node. The discount factor $\alpha \in (0,1]$ is used to trade-off between near-term cost and long-term cost.

\subsection{Intervention Model}
The goal of intervention is to reduce the risk by allocating certain amount of resources at each time step $k$ to reduce the spreading rate $\beta_{ij}$ on particular edges or increase recover rate $\delta_i$ on particular nodes. We assume that changes on $\beta_{ij}$ and $\delta_i$ are bounded and persist for all time, i.e.,
\[
\begin{split}
        & \beta_{ij}(-1)\geq \cdots \geq \beta_{ij}(k-1)\geq \beta_{ij}(k)\geq \cdots \geq \underline{\beta}_{ij}>0, \\
        & \overline{\Delta}_{i} >\overline{\delta}_i\geq \cdots \geq \delta_{i}(k) \geq \delta_{i}(k-1) \geq \cdots \geq \delta_{i}(-1) > 0,
    \end{split}
\]
for all $k\in \N$, where $\beta_{ij}(-1)=\overline{\beta}_{ij}$ and $\delta_i(-1)=\underline{\delta}_i$ are the unmodified rates. It is easy to verify that 
\begin{equation}\label{eq:rate-bound}
    \overline{A}\geq A(k-1)\geq A(k)\geq \underline{A},\quad \forall k\in \N,
\end{equation}
where the matrices $\overline{A},\underline{A} $ are given by \eqref{eq:AB} with rates $(\overline{\beta}_{ij},\underline{\delta}_i)$ and $(\underline{\beta}_{ij},\overline{\delta}_i)$, respectively. To make $\underline{A}$ nonnegative we assume $h\overline{\Delta}_i\leq 1$ for all $i\in \mathcal{V}$.

\subsection{Resource Model}
For any two matrices $A^1\geq A^2$, we defined the resource required for the transition from $A^1$ to $A^2$ as follows
\begin{equation}\label{eq:resource}
    \begin{split}
        \Gamma(A^1,A^2)=\sum_{ij}W_{ij} U_{ij}(A^1,A^2),
    \end{split}
\end{equation}
where $W_{ij}$ are the weightings on resource costed by reducing $\beta_{ij}$ and increasing $\delta_i$, respectively. And $U$ is defined as
\begin{equation}
    U_{ij}(A^1,A^2)=\begin{cases}
    \log \bigl(\beta_{ij}^1/\beta_{ij}^2\bigr), \;(i,j)\in \mathcal{E}, \\
   \log\bigl((\overline{\Delta}_i-\delta_i^{1})/(\overline{\Delta}_i-\delta_i^2)\big),\;  i=j.
    \end{cases}
\end{equation}
Given the matrix $A^1$ and resource allocation allocation $U$, one can compute the matrix $A^2$ by 
\begin{equation}
    A^2=f(A^1,U),
\end{equation}
where $A^2$ can be computed via \eqref{eq:AB} with $\beta_{ij}^2,\delta_i^2$ given by
\[
\begin{split}
        \beta_{ij}^2=\beta_{ij}^1 e^{-U_{ij}},\quad
        \delta_{i}^2=\bigl(1-e^{-U_{ii}}\bigr)\overline{\Delta}+e^{-U_{ii}}\delta_i^1.
    \end{split}
\]
The above logarithmic resource model was proposed in \cite{somers2020sparse} to encourage sparse resource allocation. One implication of the above model is that certain proportional decrease for $\beta_{ij}$ always has the same cost. Note that it is impossible for $\beta_{ij}$ to be reduced to 0 since that would take infinite resources. 

\subsection{Problem Statement}
In this paper we are interested in resource-constrained intervention of the spreading process \eqref{eq:sys} using model predictive control. Specifically, the deliverable resource at each time step is bounded by $\overline{\Gamma}>0$. The goal of intervention is to minimize the risk $R(x(k))$ given the outbreak detection $x(k)$ at time step $k$, which can be formulated as an optimal control problem (OCP):
\begin{equation}\label{eq:ocp}
    \begin{split}
        \min_{\beta,\delta} \quad &  R(x(k)) \\
        \mathrm{s.t.}\quad &   x(\ell+1|k)=[A(\ell|k)-B(x,\ell|k)]x(\ell|k), \\ 
        & x(0|k)=x(k), \\
        & \beta_{ij}({\ell-1|k})\geq \beta_{ij}({\ell|k})\geq \underline{\beta}_{ij}, \\
        &\delta_i(\ell-1|k)\leq \delta_{i}(\ell|k)\leq \overline{\delta}_i, \\
        & \Gamma(\ell|k)\leq \overline{\Gamma},\; \ell=0,\ldots,L-1
    \end{split}
\end{equation}
where $L$ is the time horizon and $\Gamma(\ell|k)$ is short for $\Gamma(A(\ell-1|k),A(\ell|k))$. Here $\beta_{ij}(-1|k)=\beta_{ij}(k-1)$ and $\delta_i(-1|k)=\delta_i(k-1)$ refer to the rates at the previous time step. For evaluating $R(x(k))$ we use $\beta_{ij}(\ell|k)=\beta_{ij}(L-1|k)$ and $\delta_{i}(\ell|k)=\delta_i(L-1|k)$ for $\ell\geq L$.  

Note that OCP \eqref{eq:ocp} is a non-convex optimization problem due to the complex  dynamics \eqref{eq:sys}, which may be intractable for real-time intervention of large-scale spreading processes. 

\section{Dynamic Programming Relaxation}\label{sec:dp}
In this section we present a convex relaxation for OCP \eqref{eq:ocp}. The key idea is to minimize an particular upper bound of the risk \eqref{eq:risk}, which can be constructed via a dynamic programming formulation. We will show that under some mild assumptions the risk bound will decrease asymptotically after certain steps. 

\subsection{Risk Bounds}
It is well-known that positive linear systems with non-negative linear costs admit linear value (cost-to-go) functions \cite{rantzer2015scalable,berman1994nonnegative}. Similarly, we will use a dynamic-programming-alike approach to construct the risk bounds of nonlinear positive system \eqref{eq:sys}. 

\begin{proposition}\label{prop:1}
    Suppose that the rate sequences $\beta_{ij}$ and $\delta_i$ satisfy \eqref{eq:rate-bound}. If there exists a sequence of nonnegative row vectors $0\leq p(\ell|k)\in \R^{1\times n}$ such that
\begin{equation}\label{eq:cond-p-k}
    \begin{split}
        p(\ell|k)&\geq C+\alpha  p(\ell+1|k)A(\ell|k),\; \forall 0\leq \ell<L, \\
    p(L|k)&=p(L-1|k),
    \end{split}
\end{equation}
then we have $R(x(k))\leq \bar{R}(x(k)):=p(0|k)x(k) $.
\end{proposition}
The above proposition appears in \cite{somers2022multi} and its continuous-time version can be found in \cite{somers2019priority}. The underlying idea is that $\bar{R}(\hat x,\ell|k):=p(\ell|k) \hat x(\ell|k)$ provides a risk upper bound of a linear spreading process
\begin{equation}\label{eq:sys-lin}
    \hat x(\ell+1|k)=A(\ell|k)\hat x(\ell|k), \quad \hat{x}(0|k)=x(k),
\end{equation}
where the solution $\hat{x}$ bounds the state trajectory $x$ of nonlinear system \eqref{eq:sys} with the same initial condition, i.e. $\hat{x}(0|k)=x(k)\Rightarrow \hat x(\ell|k)\geq x(k+\ell)$ for all $\ell$. Thus, $\bar R(x(k)):=\bar{R}(\hat x,0|k)$ also serves as a risk bound of \eqref{eq:sys}. Since the row vector $p$ provides importance measurement of the outbreak at node $i$, we call it a \emph{priority vector}.
\subsection{Relaxed OCP}
We now present a relaxed version of OCP \eqref{eq:ocp} by minimizing the risk bounds over priority vectors, i.e.,
\begin{equation}\label{eq:dp}
    \begin{split}
        \min_{p,\beta,\delta}\quad & \bar{R}_{\epsilon} (p,x):=p(0|k)x(k)+\epsilon \sum_i p_i(0|k)\\
        \mathrm{s.t.}\quad & p(\ell|k) \geq 0,\; p(\ell|k)\geq C+\alpha p(\ell+1|k)A(\ell|k), \\
        & \beta_{ij}({\ell-1|k})\geq \beta_{ij}({\ell|k})\geq \underline{\beta}_{ij}, \\
        &\delta_i(\ell-1|k)\leq \delta_{i}(\ell|k)\leq \overline{\delta}_i, \\
        & \Gamma(\ell|k)\leq \overline{\Gamma},\; \ell=0,\ldots,L-1
    \end{split}
\end{equation}
where $\epsilon>0$ is a small constant. We will give a convex reformulation of the above problem in Section~\ref{sec:exp-cone}. 

Letting $(p^*, \beta^*, \delta^*)$ be a local optimal solution of \eqref{eq:dp}, we use $A^*$ to denote the matrix obtained by substituting $\beta^*,\delta^*$ into \eqref{eq:AB}. Since \eqref{eq:dp} can be reformulated as a convex problem, $(p^*, \beta^*, \delta^*)$ is also a global minimum. We describe the properties of optimal solutions via the following propositions. In particular, Proposition~\ref{prop:2} shows that the optimal priority vectors satisfy the Bellman's equation. Proposition~\ref{prop:3} indicates that both $A^*$ and $p^*$ are monotonically decreasing before hitting the lower bound. Proposition~\ref{prop:4} implies that the proposed controller \eqref{eq:dp} is an ``all-in'' policy, i.e., distributing all available resources at each time step.

\begin{proposition}\label{prop:2}
    The optimal solution of \eqref{eq:dp} satisfies 
    \begin{equation}\label{eq:p-star}
        p^*(\ell|k)=C+\alpha p^*(\ell+1|k)A^*(\ell|k)
    \end{equation}
    for $\ell=0,\ldots,L-1$, where $p^*(L|k)=p^*(L-1|k)$.
\end{proposition}

\begin{proposition}\label{prop:3}
    If $A(k-1)\gneqq \underline{A}$, then we have
    \begin{equation}
        A^*(\ell-1|k)\gneqq A^*(\ell|k),\quad p^*(\ell-1|k)\gneqq p^*(\ell|k)
    \end{equation}
    for $\ell=0,\ldots,\underline{L}-1$, where $\underline{L}\leq L$ is the largest integer such that $A^*(\underline{L}-1|k)\gneqq \underline{A}$. 
\end{proposition}

\begin{proposition}\label{prop:4}
    If $\Gamma(A(k-1),\underline{A})>\overline{\Gamma}$, then we have 
    \begin{equation}
        \Gamma(A(k-1),A^*(0|k))=\overline{\Gamma}.
    \end{equation}
\end{proposition}

We now investigate the closed-loop behavior under the MPC controller \eqref{eq:dp}, i.e., implementing $\beta^*(0|k),\delta^*(0|k)$ at time step $k$ and resolving \eqref{eq:dp} when $x(k+1)$ is available. We first make the following assumptions.

\begin{assumption}\label{asmp:1}
The discount rate $\alpha$ satisfies $\alpha \rho(\overline{A}) < 1$ where $\rho(\cdot)$ denotes the matrix spectrum radius.
\end{assumption}

\begin{remark}
Assumption~\ref{asmp:1} implies that the risk \eqref{eq:risk} is well-defined (i.e. finite) for the linear spreading process \eqref{eq:sys-lin}, despite that the state $\hat{x}(\ell|k)$ might be unbounded. 

\end{remark}

\begin{assumption}\label{asmp:2}
The following set is nonempty:
\begin{equation}\label{eq:CA}
    \mathcal{A}:=\{A\mid \underline{A}\leq A\leq\overline{A}, \;\sum_{i}c_i\beta_{ij}<c_j\delta_j, \;\forall j\in \mathcal{V}\}.
\end{equation}
\end{assumption}

\begin{remark}
Note that $\mathcal{A}$ is convex and satisfies $A(k)\in \mathcal{A}\Rightarrow A(k+1)\in \mathcal{A}$ due to \eqref{eq:rate-bound}. Later on we will prove that for any $A\in \mathcal{A}$ the following inequality holds:
\begin{equation}\label{eq:stability}
    C-(1-\alpha)p>0,
\end{equation}
where $p=C(I-\alpha A)^{-1}$, which will be used to show the asymptotic stability of closed-loop system.
\end{remark}

\begin{theorem}\label{thm:1}
If Assumption~\ref{asmp:1} holds, then the proposed OCP \eqref{eq:dp} is recursively feasible. Furthermore, if Assumption~\ref{asmp:2} holds, then  $\bar{R}(x(k))>\bar{R}(x(k+1))$ for all $x(k)\gneqq 0$ and $k\geq K$, where $K=\left\lceil \frac{\Gamma_M}{\overline{\Gamma}} \right\rceil$ with $\Gamma_M=\liminf_{A\in\mathcal{A}} \Gamma(\overline{A},A)$.

\end{theorem}

\begin{remark}\label{rem:1}
Since $\Gamma(\overline{A},A)$ is a convex function and $\mathcal{A}$ is a convex bounded set, then $\Gamma_M$ is well-defined and serves as an estimate of the maximum resource required for making the risk bound asymptotically decreasing. We will show that $\Gamma_M$ can also be computed via exponential cone programming.
\end{remark}

\subsection{Exponential Cone Programming}\label{sec:exp-cone}

In this section we show that the proposed dynamic resource allocation problem \eqref{eq:dp} can be reformulated as a convex program, in particular exponential cone program (ECP). We refer the readers to \cite{boyd2007tutorial} for details on exponential cone programming.

First, we introduce new decision variables as follows
\begin{equation}\label{eq:variable}
    y_i^{\ell}=\log p_i(\ell),\quad  
    U_{ij}^\ell=U(A(\ell-1),A(\ell))
\end{equation}
with $\ell=0,\ldots, L-1$, where the time step $k$ is omitted for simplicity. The following convex programming is an equivalent formulation of \eqref{eq:dp}:
\begin{subequations}\label{eq:ecp_formulation}
    \begin{align}
        \min_{r,y,U}\; & r \nonumber\\
        \mathrm{s.t.}\;  
        & \log\left(\sum_i \exp \bigl(\log(x_i+\epsilon)+y_i^0-r \bigr)\right)\leq 0 \label{eq:cons-cost} \\
        &\log \bigl( \sum_{i} \exp \bigl( y_i^{\ell+1}-y_j^{\ell}+\log a_{ij}-\sum_{s=0}^{\ell} U_{ij}^{s}\bigr) \nonumber\\
        &\quad+\exp\bigl( y_j^{\ell+1}-y_j^{\ell}+\log b_j-\sum_{s=0}^\ell U_{jj}^s\bigr) \nonumber\\
        &\quad +\exp(\log c_j - y_j^\ell) \nonumber\\
        &\quad+\exp\left(y_j^{\ell+1}-y_j^{\ell}+\log d_j\right)\bigr)
        \leq 0,\label{eq:cons-pA}\\
        &\quad\quad j=1,\ldots,n,\; \ell=0,\ldots,L-1\nonumber \\
        & U^\ell\geq 0,\; \sum_{\ell=0}^{L-1}U^l\leq U(A(-1),\underline{A}), \label{eq:cons-res}\\
        &\sum_{ij} W_{ij}U_{ij}^{\ell}\leq \overline{\Gamma} ,\; \ell=0,\ldots,L-1 \nonumber 
    \end{align}
\end{subequations}
where $a_{ij}=\alpha h \beta_{ij}(k-1)$, $b_j=\alpha h(\overline{\Delta}_j-\delta_j(k-1))$, $d_j=\alpha(1-h\overline{\Delta}_j)$, and $y^L=y^{L-1}$.
\begin{remark}
Constraint \eqref{eq:cons-cost} implies $\bar R_\epsilon(p,x)\leq e^r$ while \eqref{eq:cons-pA} is equivalent to the Bellman inequality \eqref{eq:cond-p-k}, see \cite{somers2022multi} for details. The $\ell_1$-norm constraint in \eqref{eq:cons-res}, which is based on our logarithmic resource model \eqref{eq:resource}, can encourage sparse resource allocation. One can apply the reweighted $\ell_1$ optimization approach \cite{candes2008enhancing} to further minimize the number of nodes and edges with non-zero resource allocated, which is beneficial when resources cannot be distributed widely.
\end{remark} 

We also give an ECP for computing $\Gamma_M$ in Thm.~\ref{thm:1}:
\begin{subequations}
    \begin{align}
        \min_{U}\; & \sum_{ij} W_{ij} U_{ij} \nonumber\\
    \mathrm{s.t.}\; & 0\leq U\leq U(\overline{A},\underline{A}) \label{eq:Gamma-M-1}\\
    &\log\bigl ({\textstyle \sum}_i \exp (\log(c_i\overline{\beta}_{ij}/(c_j\overline{\Delta}_j))-U_{ij}) \nonumber\\
    &\;+\exp(\log((\overline{\Delta}_j-\underline{\delta}_j)/\overline{\Delta}_j)-U_{jj})+\epsilon_2\bigr)\leq 0\label{eq:Gamma-M-2}
    \end{align}
\end{subequations}
where $\epsilon_2>0$ is a small constant and $\Gamma_M$ is the optimal objective value.
The set $\mathcal{A}$ is transformed into Constraints \eqref{eq:Gamma-M-1} - \eqref{eq:Gamma-M-2} since $\eqref{eq:Gamma-M-1} \Leftrightarrow \underline{A}\leq A=f(\overline{A},U)\leq \overline{A}$ and
\[
\begin{split}
    \eqref{eq:Gamma-M-2}&\Rightarrow \sum_i\frac{c_i\overline{\beta}_{ij}e^{-U_{ij}}}{c_j\overline{\Delta}_j}+\frac{c_j(\overline{\Delta}_j-\underline{\delta}_j)e^{-U_{jj}}}{c_j\overline{\Delta}_j}<1 \\
    &\Rightarrow \sum_i\frac{c_i\beta_{ij}}{c_j\overline{\Delta}_j}+\frac{c_j(\overline{\Delta}_j-\delta_j)}{c_j\overline{\Delta}_j}<1\\
    &\Rightarrow \sum_{i}c_i\beta_{ij}<c_j\delta_j.
\end{split}
\]

\section{Proofs}\label{sec:proof}
For the sake of simplicity, we use $L_e=L-1$ to denote the last step in the time horizon and omit the dependency of time step $k$ if it does not cause any confusion. 

\subsection{Proof of Proposition~\ref{prop:1}}\label{pf:prop-1}
We first show that $\bar{R}(\hat x, \ell)=p(\ell)\hat x(\ell)$ with $p(\ell)=p(L_e)$ for all $\ell> L_e$ is a risk upper bound of the linear spreading process \eqref{eq:sys-lin}. Since $\hat x(\ell)\geq 0$, Inequality \eqref{eq:cond-p-k} implies 
\begin{equation}
    \bar{R}(\hat x, \ell)\geq C x(\ell)+\alpha \bar R (\hat{x}, \ell+1), \; \forall 0\leq \ell\leq L_e.
\end{equation}
For the case where $\ell > L_e$ we have
\begin{equation}
    \begin{split}
        \bar{R}(\hat{x},\ell) &= p(L_e)\hat{x}(\ell) \geq (C+\alpha p(L_e+1)A(L_e))\hat x(\ell) \\
        &\geq C\hat x(\ell)+\alpha p(L_e+1)A(\ell)\hat x(\ell) \\
        &= C\hat x(\ell)+\alpha p(\ell+1)A(\ell)\hat x(\ell) \\
        &= C x(\ell)+\alpha \bar R (\hat{x}, \ell+1)
    \end{split}
\end{equation}
where the last inequality follows from \eqref{eq:rate-bound}. By telescoping sum we obtain $\bar R(\hat{x},\ell)\geq \sum_{t=0}^{T-1}\alpha^t x(t+\ell)+\alpha^T \bar{R}(\hat{x},\ell+T)$ for any $T\in \N$. This further implies $\bar{R}(\hat{x},\ell)\geq R(\hat{x},\ell)$ as $T\rightarrow \infty$, where $R(\hat{x},\ell)=\sum_{t=0}^\infty \alpha^t C\hat{x}(t+\ell)$.

Secondly, we show that $R(\hat x, \ell)$ also bounds the risk of the nonlinear system \eqref{eq:sys}. By setting $\hat{x}(0)=x(0)=x$ we can have the following induction: if $\hat{x}(\ell)\geq x(\ell)$ then
\[
\begin{split}
    \hat{x}(\ell+1)&=A(\ell)\hat{x}(\ell)\geq A(\ell)x(\ell) \\
    &\geq [A(\ell)-B(x,\ell)]x(\ell)=x(\ell+1)
\end{split}
\]
where inequalities follows by nonnegative $A,B$. Now we can show that $R(x)\leq  R(\hat{x},0)\leq \bar{R}(\hat{x},0)=\bar{R}(x)$.

\subsection{Proof of Proposition~\ref{prop:2}}\label{pf:prop-2}
Suppose that $0\leq s\leq L_e$ is the smallest integer such that \eqref{eq:p-star} does not hold, i.e.,
\begin{equation}\label{eq:pf2-pstar}
    p^*(s)\gneqq C+\alpha p^*(s+1)A^*(s).
\end{equation}
If $s=0$, we construct a new sequence $\tilde{p}$ with $\tilde{p}(\ell)=p^*(\ell)$ for $\ell>0$ and $\tilde{p}(0)=C+\alpha p^*(1)A^*(0)\lneqq p^*(0)$. We have
\begin{equation}\label{eq:contradict-1}
    \begin{split}
        \bar R_\epsilon(p^*,x)-\bar{R}_\epsilon(\tilde{p},x)&=\sum_i(p_i^*(0)-\tilde{p}_i(0))(x_i+\epsilon) \\
        &\geq \epsilon \sum_i (p_i^*(0)-\tilde{p}_i(0))>0,
    \end{split}
\end{equation}
which contradicts with the optimality of $p^*$. 

For the case of $s \geq1$, we construct a sequence $\tilde{p}$ with $\tilde{p}(\ell)=p^*(\ell)$ for $s<\ell \leq L_e$ and 
\begin{equation}
    \tilde{p}(\ell)=C+\alpha \tilde{p}(\ell+1)A^*(\ell),\quad \ell=s,\ldots,0.
\end{equation}
By comparing $p^*$ and $\tilde{p}$ we have
\[
    \begin{split}
        p^*(\ell)-\tilde{p}(\ell)=\alpha [p^*(\ell+1)-\tilde{p}(\ell+1)]A(\ell),\; \forall  0\leq \ell<s.
    \end{split}
\]
Recursively applying the above relationship yields
\begin{equation}
    p^*(0)-\tilde{p}(0)=\alpha^{s}[p^*(s)-\tilde{p}(s)]\prod_{j=1}^{s}A(s-j).
\end{equation}
By the construction we have $p^*(s)\gneqq \tilde{p}(s)$. Moreover, since $A(\ell)\geq 0$ has positive diagonal, we can conclude that $p^*(0)\gneqq \tilde{p}(0)$. The rest is same as the case of $s=0$.

\subsection{Proof of Proposition~\ref{prop:3}}\label{pf:prop-3}

For the case where $\underline{L}<L$, we can treat \eqref{eq:dp} as an OCP with horizon of $\underline{L}$ as no resource can be distributed at time step $\ell\geq \underline{L}$. Therefore, we only give the proof for $\underline{L}=L$.

We first prove that $A^*(\ell-1)\gneqq A^*(\ell)$ for all $0\leq \ell \leq L_e$ by contradiction. Let $0\leq s \leq L_e$ be the largest integer such that $A^*(s-1)=A^*(s)$, i.e., there is no resource allocated at time step $s$. We will construct feasible sequences $\tilde{A}$ and $\tilde{p}$ which yield a smaller cost than $\bar{R}_\epsilon(p^*,x)$. 

We construct $\tilde{A}$ by allocating certain amount of additional resources to certain link/node at the time step $s$. Since $A^*(L_e)\gneqq \underline{A}$, without loss of generality we assume that $\beta_{ij}^*(L_e)> \underline{\beta}_{ij}$ for some link $(i,j)\in \mathcal{E}$. We now define a new sequence
\begin{equation}\label{eq:beta-tilde}
    \tilde{\beta}_{ij}(\ell)=\begin{cases}
    \beta_{ij}^*(\ell), & 0\leq \ell < s\\
    \eta \beta_{ij}^*(\ell), & s\leq \ell \leq L_e \\
    \end{cases}
\end{equation}
where $ \eta=\max\left(e^{-\overline{\Gamma}},\frac{\underline{\beta}_{ij}}{\beta_{ij}^*(L_e)}\right) $. Note that $\tilde{\beta}_{ij}$ is decreasing and satisfies $\tilde{\beta}_{ij}(L_e)\geq \underline{\beta}_{ij} $. The resource allocation does not change except the time step $s$, which satisfies
\[
\tilde{\Gamma}(s)=w_{ij}\log\frac{\beta_{ij}^*(s-1)}{\eta \beta_{ij}^*(s)}=-w_{ij}\log\eta\leq \overline{\Gamma}.
\]
Thus, $\tilde{A}$ is a feasible solution. 

We now construct a feasible $\tilde{p}$ based on $\tilde{A}$ as follows
\begin{equation}\label{eq:p-tilde}
    \tilde{p}(\ell)=C+\alpha \tilde{p}(\ell+1)\tilde{A}(\ell), \quad \ell=L_e-1,\ldots,0
\end{equation}
where $\tilde{p}(L_e)=C\bigl(I-\alpha \tilde{A}(L_e)\bigr)^{-1} $. From \eqref{eq:p-star} and \eqref{eq:p-tilde} we can obtain
\[
    \begin{split}
        p^*(& L_e)-\tilde{p}(L_e) \\
        &\quad =C(I-\alpha A^*(L_e))^{-1}-C(I-\alpha \tilde{A}(L_e))^{-1} \\
        &\quad= \alpha C(I-\alpha A^*(L_e))^{-1}[A^*(L_e)-\tilde{A}(L_e)] \\
        &\quad\quad\quad (I-\alpha \tilde{A}(L_e))^{-1} \gneqq 0.
    \end{split}
\]
Furthermore, by induction we have 
\begin{equation}\label{eq:p-1}
    \begin{split}
        p^*(\ell)-\tilde{p}(\ell)
        =&\alpha [p^*(\ell+1)A^*(\ell)-\tilde{p}(\ell+1)\tilde{A}(\ell)] \\
        =&\alpha p^*(\ell+1)(A^*(\ell)-\tilde{A}(\ell))+\\
        &\alpha (p^*(\ell+1)-\tilde{p}(\ell+1))\tilde{A}(\ell) \\
        &\gneqq \alpha (p^*(\ell+1)-\tilde{p}(\ell+1))\tilde{A}(\ell)\gneqq 0
    \end{split}
\end{equation}
where $s\leq \ell \leq L_e$. For $0\leq \ell <s $, from \eqref{eq:p-star}, \eqref{eq:beta-tilde} and \eqref{eq:p-tilde} we can obtain
\begin{equation}\label{eq:p-2}
    p^*(\ell)-\tilde{p}(\ell)=\alpha (p^*(\ell+1)-\tilde{p}(\ell+1))A^*(\ell).
\end{equation}
Recursively applying \eqref{eq:p-1} and \eqref{eq:p-2} yields
\[
\begin{split}
    p^*&(0)-\tilde{p}(0)=\\
    &\alpha^{L_e} (p^*(L_e)-\tilde{p}(L_e))\times \prod_{j=0}^{s-1} A^*(j) \times \prod_{j=s}^{L_e}\tilde{A}(j) \gneqq 0,
\end{split}
\]
which leads to a contradiction due to \eqref{eq:contradict-1}.

We prove the second part of this proposition, i.e., $p^*(\ell-1)\gneqq p^*(\ell)$ for $0\leq \ell \leq L_e$. From \eqref{eq:p-star} we have
\[
\begin{split}
    p^*(L_e&-1)-p^*(L_e) \\
    =& \alpha [p^*(L_e)A^*(L_e-1)-p^*(L_e+1)A^*(L_e)] \\
    =& \alpha p^*(L_e)[A^*(L_e-1)-A^*(L_e)]\gneqq 0
\end{split}
\]
where the second equality is due to $p^*(L_e)=p^*(L_e+1)$. Now we can obtain $p^*(\ell-1)\gneqq p^*(\ell)$ recursively since
\[
    \begin{split}
        p^*(\ell-1)-p^*(\ell)&=\alpha[p^*(\ell)A^*(\ell-1)-p^*(\ell+1)A^*(\ell)] \\
        &=\alpha[p^*(\ell)(A^*(\ell-1)-A^*(\ell))+\\
        &\quad\quad (p^*(\ell)-p^*(\ell+1))A^*(\ell)] \\
        & \gneqq \alpha (p^*(\ell)-p^*(\ell+1))A^*(\ell), \;\forall \ell<L_e.
    \end{split}
\]

\subsection{Proof of Proposition~\ref{prop:4}}\label{pf:prop-4}

First, the resource allocated at time step $\ell$ is
\begin{equation}
    \Gamma^*(\ell)=\Gamma(A^*(\ell-1),A^*(\ell))=\sum_{ij}W_{ij} U_{ij}^{*}(\ell),
\end{equation}
where $W$ is the resource weighting matrix. Supposing that $\gamma = \overline{\Gamma} -\Gamma^*(0) >0$, we aim to construct a new solution that leads to a lower risk bound. Let $\underline{L}\in[0,L] $ be the largest integer such that $A^*(\underline{L}-1)\gneqq \underline{A}$. Note that such $\underline{L}$ exists; otherwise we have $\Gamma(A^*(-1),A^*(0))=\Gamma(A^*(-1),\underline{A})>\overline{\Gamma}$. Similar to the proof of Proposition~\ref{prop:3}, we only consider the case where $\underline{L}=L$. 

First, we construct a nonnegative matrix $\widehat{U}$ with $\widehat{U}_{ij}=0$ for $i\neq j$ and $(i,j)\notin \mathcal{E}$ such that 
\[
f(A^*(L-1),\widehat{U}) \geq \overline{A} \quad \text{and}\quad \sum_{ij} W_{ij} \widehat{U}_{ij}\leq \gamma'
\]
where $\gamma'=\min\{\Gamma(A^*(L_e),\underline{A}),\gamma\}$. Then we consider a new sequence $\widetilde{U}$ with $\widetilde{U}(0)=U^{*}(0)+\widehat{U}$ and $\widetilde{U}(\ell)=U^{*}(\ell)$ for $\ell=1,\ldots,L_e$. Then, a sequence $\tilde{A}$ can be constructed via
\begin{equation}\label{eq:A-tilde}
    \tilde{A}(\ell)=f(\tilde{A}(\ell-1),\widetilde{U}(\ell))=f(A^*(\ell),\widehat{U})
\end{equation}
where $\tilde{A}(-1)=A^*(-1)$. It is easy to verify that $\tilde{A}$ is a feasible solution since
\[
\begin{split}
    \tilde\Gamma(0)&=\Gamma^*(0)+\sum_{ij} W_{ij} \widehat{U}_{ij}=\overline{\Gamma}-\gamma+\gamma'\leq \overline{\Gamma}, \\
    \tilde \Gamma(\ell)&=\Gamma^*(\ell)\leq \overline{\Gamma},\quad \ell=1,\ldots,L_e.
\end{split}
\]
By construction \eqref{eq:A-tilde} we have $A^*(\ell)\gneqq \tilde{A}(\ell)$. Thus, similar technique in the proof of Proposition~\ref{prop:3} can be used to obtain a contradiction.

\subsection{Proof of Theorem~\ref{thm:1}}\label{pf:thm-1}

\paragraph{Recursive feasibility} Assumption~\ref{asmp:1} implies that $I-\alpha A(-1)$ has a nonnegative inverse \cite[p. 113]{berman1994nonnegative}. Since $A$ is non-increasing, we have $ \rho(A(k-1))\leq \rho(A(k))$ (\cite[p. 15]{berman1994nonnegative}) and thus $(I-\alpha A(k))^{-1}\gneqq 0$ for all $k$. 

Suppose that OCP \eqref{eq:dp} is feasible at the time step $k$. We construct a candidate sequence for the time step $(k+1)$ as follows
\begin{equation}
    \begin{split}
        A(\ell|k+1)&=
    \begin{cases}
    A^*(\ell+1|k), & 0\leq \ell \leq L-2 \\
    A^*(L-1|k), & \ell=L-1
    \end{cases} \\
    p(\ell|k+1)&=
    \begin{cases}
    p^*(\ell+1|k), & 0\leq \ell \leq L-2 \\
    p^*(L-1|k), & \ell=L-1
    \end{cases}
    \end{split}
\end{equation}
i.e., one time-step shift of the previous optimal solution. Note that $p$ is a nonnegative sequence and the sequence $A$ satisfies the non-increasing constraint \eqref{eq:rate-bound} and resource budget $\Gamma(\ell|k+1)\leq \overline{\Gamma}$. Thus, $(p,A)$ is a feasible solution.

\paragraph{Stability} We first prove that Assumption~\ref{asmp:2} implies \eqref{eq:stability} since
\begin{equation}\label{eq:stable-cond-pf}
    \begin{split}
    \eqref{eq:CA}&\Rightarrow \alpha h(c_j\delta_j- {\textstyle\sum}_{i} c_i\beta_{ij})>0 \\
    &\Rightarrow \alpha C(I-A)>0 \\
    &\Rightarrow C(I-\alpha A)>(1-\alpha) C\\
    &\Rightarrow C-(1-\alpha)C(I-\alpha A)^{-1}>0\Rightarrow\eqref{eq:stability}.
\end{split}
\end{equation}

Then, we show  $A^*(0|k)\in \mathcal{A}$ for $k\geq K$. From Proposition~\ref{prop:4} we have
\[
\Gamma(\overline{A},A^*(0|k))=(k+1)\overline{\Gamma}> \liminf_{A\in\mathcal{A}} \Gamma(\overline{A},A).
\]
This also implies that Condition~\eqref{eq:stability} holds for $p^*(0|k)$ with $k\geq K$. Finally, the risk bounds are asymptotically decreasing for $k\geq K$ and $x(k)\gneqq 0$ since
\begin{equation}\label{eq:lyapunov-ineq}
    \begin{split}
        \bar{R}(x(k))&-\bar{R}(x(k+1)) \\
        &=p^*(0|k)x(k)-p^*(0|k+1)x(k+1) \\
        &\geq p^*(0|k)x(k)-p^*(1|k)x(k+1) \\
        &\geq (p^*(0|k)-p^*(1|k)A^*(0|k))x(k) \\
        &= \frac{1}{\alpha}[C-(1-\alpha)p^*(0|k)]x(k) >0.
    \end{split}
\end{equation}

\section{An Illustrative Example}\label{sec:example}
\begin{figure}[!ht]
        \begin{center}
            \includegraphics[width = 0.48\textwidth]{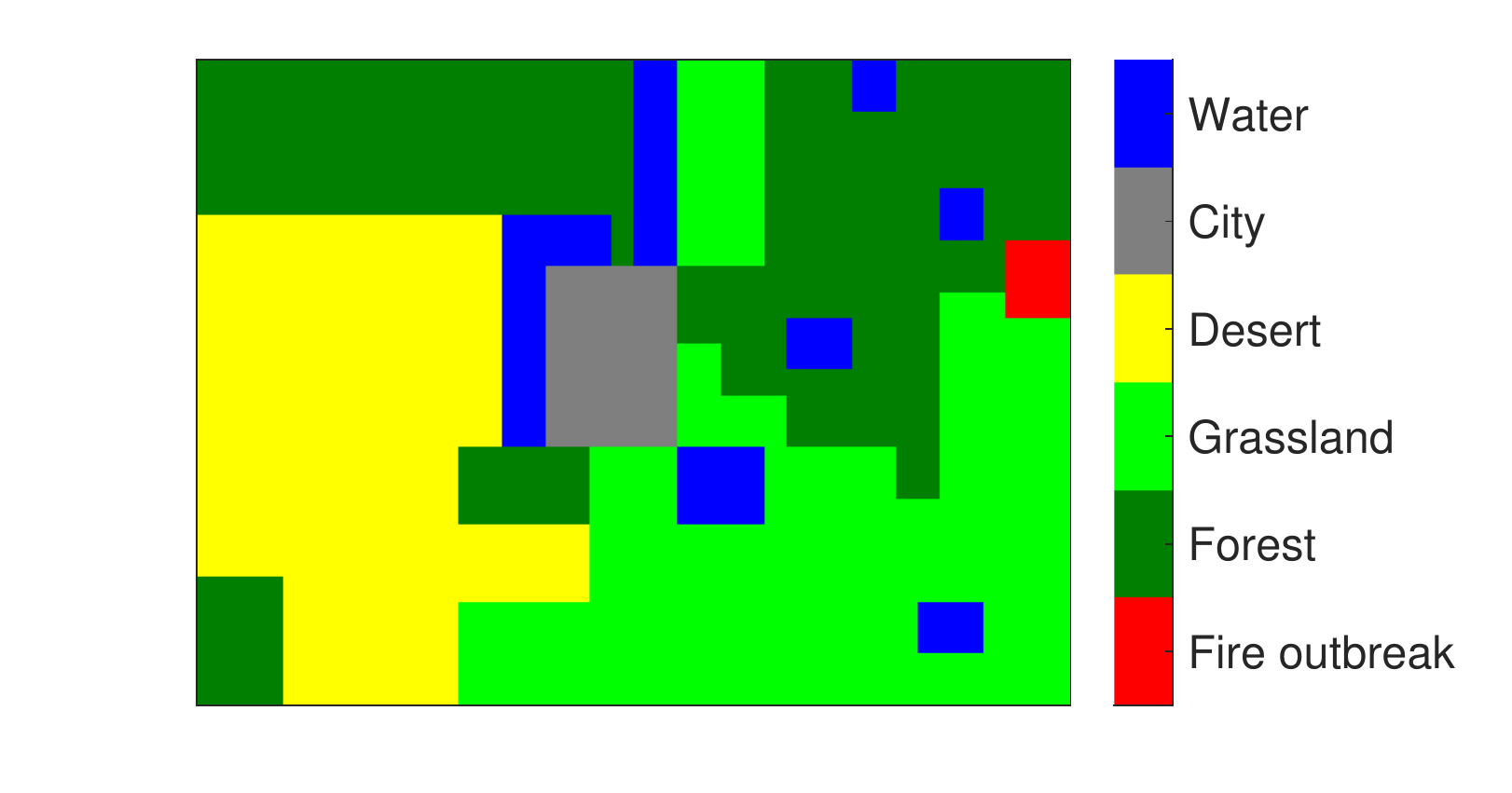} \caption{Fictional landscape with different type of areas and the location of fire outbreak.}
            \label{fig:vegetation_with_fire_outbreak}
        \end{center}
    \end{figure}

    Let us consider an example of a wildfire which can spread on a fictional landscape with different area types, as shown in Fig.~\ref{fig:vegetation_with_fire_outbreak}. We represent this landscape as a graph $\mathcal{G}(\mathcal{V},\mathcal{E})$ with a total number of nodes $n = 1000$. We consider that the fire can spread in vertical, horizontal, and diagonal direction i.e., the edge set $\mathcal{E}$ connects each node to a maximum of 8 neighbouring nodes. We also consider the precise location of a fire outbreak as shown in Fig.~\ref{fig:vegetation_with_fire_outbreak}. For different vegetation type, the rates of spreading are determined according to the wildfire models in \cite{karafyllidis1997model} and \cite{alexandridis2008cellular}. The spreading rate of an edge is given by $\beta = \beta_b \beta_{veg} \beta_{w}$, where $\beta_b = 0.5$ is the baseline, and $\beta_{veg} = 0.1,~1,~1.4$ for desert, grassland and eucalyptus forest, respectively. Water is considered unburnable, hence those edges are removed. Following \cite{alexandridis2008cellular}$, \beta_{w}$ is calculated for a northeasterly wind with a wind speed of $4m/s$. The spreading rate is adjusted, following \cite{karafyllidis1997model}, for a spread between the diagonally connected nodes. For all nodes $i\in\mathcal{V}$, the recovery rate is set as $\delta = 0.5$.
    
    We implement our proposed MPC approach to minimize the risk associated with the outbreak by allocating budgeted resources at each time step. We set the weightings $w_{ij}~=~1,~ \forall~(i,j)\in\mathcal{E}$ which indicates an equal cost to apply resources to any edge. 
    We assign a higher cost of $c_j = 1$ for all the city nodes, which could reflect either a higher risk to human life or a higher economic cost of fire reaching the city. The cost for the rest of the nodes is set as $c_j = 0.001$. The discount rate is set as $\alpha~=~1/(0.05+\rho(A))$ such that Assumption~\ref{asmp:1} holds.
    
\begin{figure}[!ht]
        \begin{center}
            \includegraphics[width = 0.48\textwidth]{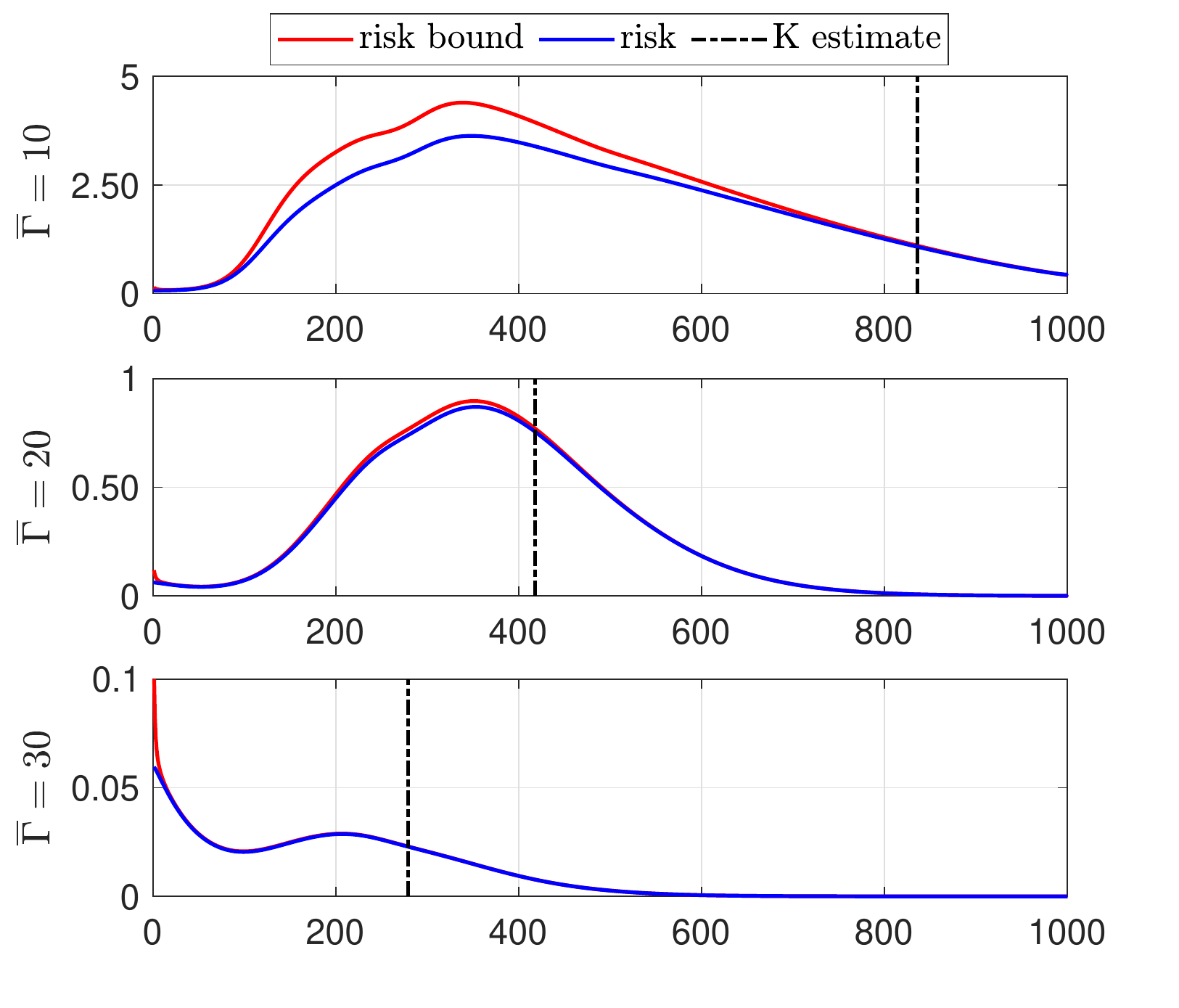} \caption{Evolution of risk $R(x(k))$, the upper bound on risk $\bar{R}(p,x)$ and the value of $K-$estimate for different values of resource bound $\overline{\Gamma}$.}
            \label{fig:Horizon_2}
        \end{center}
    \end{figure}

    The performance of our MPC method with $L=1$ and different resource bounds $\overline{\Gamma} = 10,~20,$ and $30$, is shown in Fig.~\ref{fig:Horizon_2}. 
    It shows the evolution of risk $R(x(k))$ given in \eqref{eq:risk} as well as the risk bound $\bar{R}(\hat x,0|k)$ (see Proposition~\ref{prop:1}), as the number of time steps increases.
    It can be seen that for all three cases, our MPC method is able to reduce the risk associated with the outbreak as time increases. However, for a smaller resource budget of $\overline{\Gamma} = 10$, the risk tends to increase first which reflects the unavailability of sufficient resources to suppress the risk at initial time steps. On the other hand, with $\overline{\Gamma} = 30$, the risk starts to decrease almost instantly. Note that in all three case, the value of risk $R(x(k))$ is always less than the risk bound $\bar{R}(\hat x,0|k)$. It also shows the estimated value of $K$ such that for $k\geq K$, the risk bound decreases asymptotically (see Theorem~\ref{thm:1} and Remark~\ref{rem:1}). Note that, for instance $\overline{\Gamma} = 10$, the $K-$estimate is $836$ but the value of risk and the risk bound start to decrease asymptotically after $k=400$.
    
\begin{figure}[!ht]
        \begin{center}
            \includegraphics[width = 0.48\textwidth]{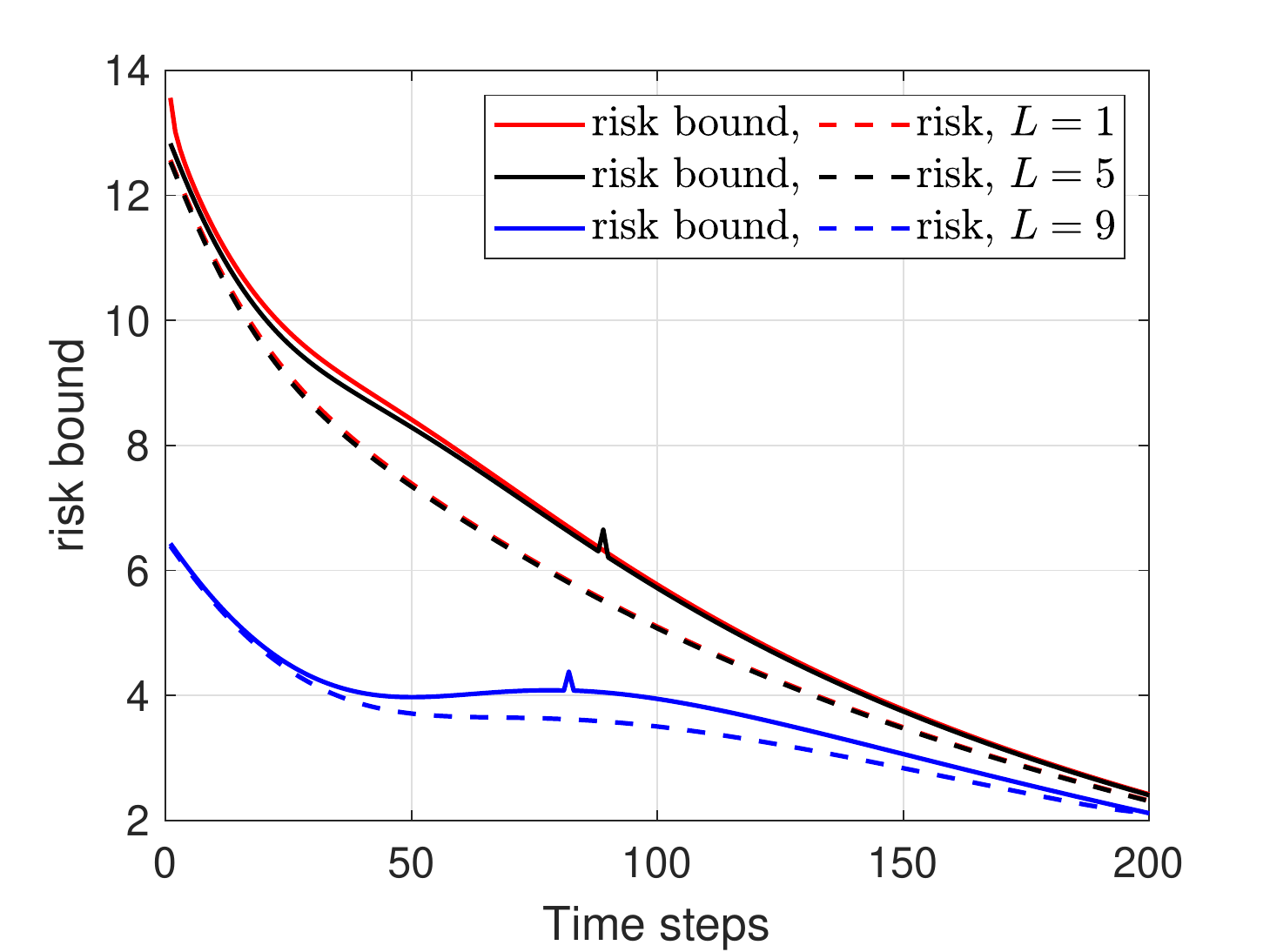} \caption{Effect of increasing $L$ on risk $R(x(k))$ and the risk bound $\bar{R}(p,x)$, with $\overline{\Gamma} = 40$ and $25\%$ of nodes being infected at time $k=0$.}
            \label{fig:risk_vs_L}
        \end{center}
    \end{figure}

    Fig.~\ref{fig:risk_vs_L} illustrates the effect of increasing $L$ on the evolution of the risk $R(x(k))$ and the upper bound on risk $\bar{R}(p,x)$ while we set the resource budget as $\overline{\Gamma} = 40$ and we consider that the fire outbreak affects $25\%$ of nodes at time $k=0$. It indicates that a larger value of $L$ performs better by resulting in a smaller value of risk at the initial time steps. 

\bibliographystyle{IEEEtran}
\bibliography{ref}

\end{document}